\documentclass[11pt,a4paper]{article}
\usepackage[margin=27mm]{geometry}
\usepackage[T1]{fontenc}
\usepackage{lmodern}
\usepackage{amsmath,amssymb,amsthm}
\IfFileExists{microtype.sty}{\usepackage{microtype}}{}
\usepackage[hidelinks]{hyperref}
\hypersetup{pdftitle={Extremal graphs for Alon's graph coupon-collector conjecture},
 pdfauthor={Xingchao Deng, Guoqing Wang}}
\newtheorem{theorem}{Theorem}[section]
\newtheorem{proposition}[theorem]{Proposition}
\newtheorem{lemma}[theorem]{Lemma}
\newtheorem{corollary}[theorem]{Corollary}
\newtheorem*{conjecture}{Conjecture}
\theoremstyle{remark}
\newtheorem{remark}[theorem]{Remark}
\newcommand{\cc}{\operatorname{cc}}
\newcommand{\E}{\mathbb{E}}
\newcommand{\Pp}{\mathbb{P}}
\newcommand{\Exp}{\operatorname{Exp}}
\newcommand{\dd}{\,\mathrm{d}}

\begin{document}
\title{\textbf{Extremal graphs for Alon's graph coupon-collector conjecture}}
\author{\small\textsc{Xingchao Deng and Guoqing Wang}}
\date{}
\maketitle

\begin{abstract}
Let \(G\) be a finite simple graph without isolated vertices.
Starting with every vertex uncovered, at each step sample a vertex
uniformly and independently and cover one of its uncovered neighbors,
chosen uniformly, if any. Let \(\cc(G)\) denote the expected number
of steps required to cover all vertices. Alon conjectured that
disjoint unions of \(K_{d+1}\) maximize \(\cc(G)\) among
\(d\)-regular graphs of a prescribed order whenever such unions exist.
In this paper, we disprove this conjecture for every fixed \(d\ge2\)
and give a complete characterization of the extremal graphs for
\(d=2\) at every order \(n\ge3\). In particular, all but at most
six vertices of any such extremal graph lie in four-cycle components.
\end{abstract}

\section{Introduction}\label{sec:introduction}

The classical coupon-collector problem asks how many independent
random samples are needed to collect every coupon type.
Alon~\cite{Alon2008} considered the following graph version.
On a finite simple graph \(G\) without isolated vertices, all vertices
are initially uncovered. At each step, a vertex is sampled uniformly
and independently of the previous samples. If it has uncovered
neighbors, one of them is chosen uniformly and covered; otherwise
the covered set remains unchanged. The sampled vertex may itself be
covered or uncovered. Let \(\cc(G)\) denote the expected number of
steps until all vertices are covered, and write \(mF\) for the
disjoint union of \(m\) copies of \(F\). Alon conjectured that, among
regular graphs of prescribed order and degree, this expectation is
maximized by disjoint unions of cliques.

\begin{conjecture}[Alon {\cite[Section~6, Remark~3]{Alon2008}}]
Let \(d,n\) be positive integers such that \(d+1\) divides \(n\).
For every finite simple \(d\)-regular graph \(G\) on \(n\) vertices,
\[
 \cc(G)\leq\cc\!\left(\frac{n}{d+1}K_{d+1}\right).
\]
\end{conjecture}

Adler, Halperin, Karp, and Vazirani
\cite[Theorem~4.2]{Adler2003} established linear cover time for the
hypercube, while Alon \cite[Theorems~6.1 and~6.2]{Alon2008} proved a
general lower bound and upper bounds for expanders. Dimitrov and
Plaxton \cite{DimitrovPlaxton2013} obtained cover-time bounds of the
correct order for related processes in which neighbors are selected
using random vertex ranks. For further background on coupon collection,
see the survey
of Boneh and Hofri \cite{BonehHofri1997} and the references therein.

In this paper, we disprove Alon's conjecture in every fixed degree
at least two. In degree two, we determine all graphs maximizing
\(\cc(G)\) at each prescribed order; the smallest counterexample
has nine vertices. The complete extremal characterization is given
by the following theorem.

\begin{theorem}\label{thm:main}
For \(n\geq3\), define
\[
 G_n^\star=
 \begin{cases}
  \dfrac n4 C_4,&n\equiv0\pmod4,\\[4pt]
  C_5\sqcup\dfrac{n-5}{4}C_4,&n\equiv1\pmod4,\\[4pt]
  C_6\sqcup\dfrac{n-6}{4}C_4,&n\equiv2\pmod4,\\[4pt]
  C_3\sqcup\dfrac{n-3}{4}C_4,&n\equiv3\pmod4,
 \end{cases}
\]
where \(\sqcup\) denotes disjoint union. For every finite simple
\(2\)-regular graph \(G\) on \(n\) vertices,
\begin{equation}\label{eq:extremal-expectation}
 \cc(G)\leq\cc(G_n^\star).
\end{equation}
Equality holds precisely when \(G\cong G_n^\star\), or when
\(n\equiv2\pmod4\) and
\[
 G\cong 2C_3\sqcup\frac{n-6}{4}C_4.
\]
Moreover, let \(T_H\) be the completion time of the same covering
process on a graph \(H\), with each vertex sampled at the rings of
an independent rate-\(1\) Poisson clock. Then
\begin{equation}\label{eq:extremal-cdf}
 \Pp(T_G\leq t)\geq\Pp(T_{G_n^\star}\leq t)
 \qquad(t>0),
\end{equation}
with equality at any positive \(t\) for exactly the same graphs.
\end{theorem}

For fixed degree \(d\ge2\) and a fixed bound of at least \(2d\)
on component order, unions of \(K_{d,d}\) uniquely maximize
\(\cc(G)\) at all sufficiently large orders divisible by \(2d\).
The results for general degree also identify the neighborhood structure
governing the completion-time tail and disprove the uniform additive
estimate suggested alongside Alon's conjecture.

\section{Preliminaries and component distributions}\label{sec:poisson}

Let \(N_G\) be the number of discrete steps required to cover
\(G\), so that \(\cc(G)=\E N_G\). This uniform neighbor-selection
rule is called process UNI in \cite{DimitrovPlaxton2013}.
Give each vertex an independent Poisson clock of rate \(1\), and apply
the same covering rule whenever its clock rings. Let \(T_G\) denote the
completion time of this continuous-time process. We write
\(\Exp(a)\) for an exponential random variable of rate \(a\).
All logarithms are natural. We use standard facts about Poisson clocks
and finite-state continuous-time Markov chains; see
\cite[Sections~2.1, 2.2, 2.4, 4.1, and~4.2]{Durrett2016}.
The tail-integral and exponential-distribution identities below can be
found in \cite[Section~2.1 and Appendix~A.3]{Durrett2016}.

\begin{lemma}
\label{lem:poisson-clocks}
Let independent Poisson clocks have rates
\(\lambda_1,\ldots,\lambda_s>0\), and put
\(\Lambda=\sum_{i=1}^s\lambda_i\). Their superposition has rate
\(\Lambda\). Its labels are independent, with probabilities
\(\lambda_i/\Lambda\), and independent of the interarrival times.
After an almost surely finite stopping time, the clocks restart
independently with the same rates, independently of the past.
\end{lemma}

\begin{lemma}
\label{lem:finite-state-chains}
Finite-state continuous-time Markov chains with the same generator
and initial distribution have the same law. At a nonabsorbing state,
the holding time is exponential with the total transition rate and
independent of the next state.
\end{lemma}

\begin{lemma}
\label{lem:tail-exponential}
For a nonnegative random variable \(T\),
\begin{equation}\label{eq:tail}
 \E T=\int_0^\infty \Pp(T>t)\dd t.
\end{equation}
If \(X_1,\ldots,X_k\) are independent \(\Exp(\lambda)\) variables,
where \(k\ge1\) and \(\lambda>0\), then
\[
 \Pp(X_1+\cdots+X_k>t)
 =e^{-\lambda t}\sum_{j=0}^{k-1}\frac{(\lambda t)^j}{j!}
 \qquad(t\ge0),
\]
and
\[
 \E\max_{1\le i\le k}X_i=\frac{H_k}{\lambda},
 \qquad H_k=\sum_{j=1}^k\frac1j.
\]
\end{lemma}
The last identity follows from the memoryless property: the successive
order-statistic spacings have rates
\(k\lambda,(k-1)\lambda,\ldots,\lambda\).

\begin{lemma}\label{lem:poisson}
If \(G\) has \(n\) vertices and no isolated vertices, then
\begin{equation}\label{eq:poisson}
 \cc(G)=n\E T_G.
\end{equation}
The continuous-time processes on distinct components are independent.
\end{lemma}

\begin{proof}
By Lemma~\ref{lem:poisson-clocks}, superposing the \(n\) independent
rate-\(1\) clocks gives a Poisson process of rate \(n\).
Its successive vertex labels are independent
and uniform on \(V(G)\), and independent of the interarrival times
\(Z_1,Z_2,\ldots\), which have distribution \(\Exp(n)\). Generate
the neighbor choices using further independent randomness. The
covering states at ring times then follow the discrete process, and
\[
 T_G=\sum_{j=1}^{N_G}Z_j.
\]
Since \(N_G\) depends only on the labels and neighbor choices,
conditioning on \(N_G\) gives \(\E T_G=\E N_G/n\).

The expectations are finite. From any nonempty uncovered set, choose
a vertex \(u\) and a neighbor \(v\). Sampling \(v\) produces a new
covering, so the expected number of steps until the next covering is
at most \(n\). There are exactly \(n\) coverings, hence
\(\E N_G\leq n^2\). Finally, distinct components use disjoint sets of
clocks and independent neighbor choices. Their continuous-time
processes are therefore independent.
\end{proof}

\begin{lemma}\label{lem:component-distributions}
For \(t\geq0\), define
\begin{align}
 q(t)&=1-(1+2t)e^{-2t},\label{eq:q}\\
 r(t)&=1-9e^{-2t}+(8+6t)e^{-3t}.\label{eq:r}
\end{align}
Then \(\Pp(T_{C_4}\leq t)=q(t)^2\) and
\(\Pp(T_{K_3}\leq t)=r(t)\).
Consequently, for every integer \(k\geq1\),
\begin{align}
 \Pp(T_{3kC_4}\leq t)&=q(t)^{6k},\label{eq:cycles-cdf}\\
 \Pp(T_{4kK_3}\leq t)&=r(t)^{4k}.\label{eq:triangles-cdf}
\end{align}
\end{lemma}

\begin{proof}
View \(C_4\) as \(K_{2,2}\), with vertex classes \(L\) and \(R\).
Let \(A_t\) and \(B_t\) count the clock rings in \(L\) and \(R\),
respectively, up to time \(t\). These are independent Poisson random
variables with mean \(2t\). Each ring in \(L\) covers a new vertex
of \(R\) until both vertices of \(R\) are covered, and the same
statement holds with \(L\) and \(R\) interchanged. Thus the four-cycle
is completely covered by time \(t\) if and only if \(A_t\geq2\) and
\(B_t\geq2\). It follows that
\[
 \Pp(T_{C_4}\leq t)
 =\Pp(A_t\geq2)\Pp(B_t\geq2)
 =\bigl(1-(1+2t)e^{-2t}\bigr)^2=q(t)^2.
\]

In \(K_3\), every vertex has an uncovered neighbor while at least two
vertices remain uncovered. Hence every ring is successful until the
second covering, and the first two waiting times are independent
\(\Exp(3)\) variables. After the second covering, exactly one vertex
remains uncovered. It is covered at the next ring from either of its
two neighbors, whose combined clock rate is \(2\). By the strong
Markov property in Lemma~\ref{lem:poisson-clocks}, the conditional
distribution of this last waiting
time, given the history up to the second covering, is \(\Exp(2)\).
Since this distribution does not depend on the history, the last
waiting time is independent of the first two. Let \(X\) be the sum
of the first two waiting times and \(Y\) the
last waiting time. Then \(X\) has density \(9s e^{-3s}\),
\(Y\sim\Exp(2)\), and \(X\) and \(Y\) are independent. Therefore,
\begin{align*}
 \Pp(T_{K_3}\leq t)
 &=\int_0^t 9s e^{-3s}\bigl(1-e^{-2(t-s)}\bigr)\dd s\\
 &=1-(1+3t)e^{-3t}
   -9e^{-2t}\bigl(1-(1+t)e^{-t}\bigr)\\
 &=1-9e^{-2t}+(8+6t)e^{-3t}=r(t).
\end{align*}

A disjoint union is completely covered by time \(t\) if and only if
all its components are completely covered by that time. By
Lemma~\ref{lem:poisson}, the component completion times are independent.
Thus, for \(m\) copies of a graph \(F\),
\[
 \Pp(T_{mF}\leq t)=\bigl(\Pp(T_F\leq t)\bigr)^m.
\]
Taking \(F=C_4\), \(m=3k\), and then \(F=K_3\), \(m=4k\), gives
\eqref{eq:cycles-cdf} and \eqref{eq:triangles-cdf}, respectively.
\end{proof}

For integers \(m,d\geq1\), let \(D_{m,d}\) be the number of independent
uniform draws from \(m\) types required to collect at least \(d\)
copies of each type. This is the classical repeated coupon-collector
problem.

\begin{theorem}[Newman--Shepp {\cite[Theorem~2]{NewmanShepp1960}}]
\label{thm:newman-shepp}
For each fixed \(d\geq2\),
\begin{equation}\label{eq:newman-shepp}
 \E D_{m,d}
 =m\bigl(\log m+(d-1)\log \log  m+O_d(1)\bigr)
 \qquad(m\to\infty).
\end{equation}
\end{theorem}
Erd\H{o}s and R\'enyi \cite{ErdosRenyi1961} obtained the corresponding
limiting distribution.

\section{Extremal graphs of degree two}\label{sec:extremal}

For a \(2\)-regular graph \(G\), form an auxiliary multigraph on the same
vertex set by placing one edge between the two neighbors of each
vertex of \(G\). If \(u\) is uncovered, its covering rate is
\[
 \sum_{v\sim u}\frac1{\#\{\text{uncovered neighbors of }v\}}
 =1+\frac12\#\{\text{covered auxiliary neighbors of }u\},
\]
with multiplicities counted. Thus the covering process has the same generator
as an infection process in which each vertex becomes infected
spontaneously at rate \(1\), and each infected neighbor transmits
infection at rate \(1/2\).
Both processes start from the empty set, so
Lemma~\ref{lem:finite-state-chains} shows that their evolving covered
and infected sets have the same law, with the same time parameter.
This auxiliary infection process is a special case of the discrete
Bass model. Bounds on adoption probabilities and decompositions of
susceptibility events are given in
\cite{FibichLevin2024,FibichLevinSchochet2025}.

Let \(F_m(t)\) be the probability that this infection process covers an
\(m\)-cycle by time \(t\). For \(m=2\), use two parallel edges; the two
successive waiting times then have distribution \(\Exp(2)\), so
\(F_2=q\). Set \(F_1(t)=1-e^{-t}\) for the generating functions below.
The auxiliary multigraph of an odd cycle is isomorphic to that
cycle, whereas the auxiliary multigraph of \(C_{2m}\) consists of two
disjoint \(m\)-cycles. The equality in law therefore gives
\begin{equation}\label{eq:auxiliary-cycle-law}
 \Pp(T_{C_\ell}\le t)=
 \begin{cases}
 F_\ell(t),&\ell\text{ odd},\\
 F_{\ell/2}(t)^2,&\ell\text{ even}.
 \end{cases}
\end{equation}
In particular, \(F_3=r\), and
\(T_{C_6}\stackrel{\mathrm d}=T_{2C_3}\).

\begin{lemma}\label{lem:cycle-generating-function}
For \(t>0\), put \(a=e^{-t}\). There are positive coefficients \(c_j(t)\)
such that
\begin{equation}\label{eq:cycle-generating-function}
 \sum_{m\ge1}\frac{F_m(t)}m x^m
 =-\log(1-C_t(x)),\qquad C_t(x)=\sum_{j\ge1}c_j(t)x^j,
\end{equation}
where \(c_1(t)=1-a\) and
\begin{equation}\label{eq:cycle-coefficients}
 c_j(t)=a^2\int_0^{1-a}
       \frac{v^{j-2}}{(j-2)!}\frac{1-a-v}{(a+v)^2}\dd v
 \qquad(j\ge2).
\end{equation}
All generating functions are formal power series in \(x\).
\end{lemma}

\begin{proof}
Consider an \(m\)-cycle with vertex set \(V\), where
\(m\ge2\), with two parallel edges when \(m=2\).
Start with every vertex uninfected. Place independent rate-\(1/2\)
Poisson transmission arrows on each oriented edge and independent
rate-\(1\) spontaneous-infection marks at each vertex. An arrow from
an infected vertex infects its endpoint, and infection is permanent.
Fix \(t\ge0\) and a nonempty set \(S\) of \(r\) vertices. Choose a
clockwise orientation of the cycle, and let \(k_1,\ldots,k_r\) be the
numbers of vertices outside \(S\) strictly between consecutive
vertices of \(S\). Thus
\[
 k_i\ge0,\qquad k_1+\cdots+k_r=m-r.
\]
Trace potential infection paths backward from \(S\) at time \(t\).
For \(0\le s\le t\), let \(B_s\) be the set of vertices at time
\(t-s\) from which a path, using vertical time segments and
transmission arrows, can reach \(S\) at time \(t\). Then \(B_0=S\).
On encountering an arrow \(x\to y\) during the backward exploration,
add \(x\) to \(B_s\) if \(y\) is already an ancestor. Since there is
no recovery, every vertex that enters the ancestor set remains in
\(B_s\) as \(s\) increases.

Within a gap containing \(k\ge1\) vertices, the ancestor set advances
from its two ends. Until the gap is filled, exactly two oriented
boundary edges can add another ancestor, each at rate \(1/2\).
When only one vertex remains, the two arrows capable of adding it are
still distinct. The total rate of acquiring another gap vertex is
therefore \(1\) at every stage. By the memoryless property of the
Poisson clocks, the number acquired by backward time \(s\) has the
same process law as \(\min\{P_s,k\}\), where \(P\) is a rate-\(1\)
Poisson process. Distinct gaps use disjoint families of transmission
arrows, so their growth processes are independent. This also covers
\(r=1\); for \(m=2\), the two parallel edges provide the two distinct
boundary channels. Let \(X_i(s)\) count the ancestors acquired in the
\(i\)-th gap.
The \(r\) vertices of \(S\) remain ancestors throughout, so
\[
 |B_s|=r+\sum_{i=1}^r X_i(s).
\]
Conditional on all transmission arrows, every vertex of \(S\) is
uninfected at time \(t\) exactly when the backward ancestor region
contains no spontaneous-infection mark. The marks have rate \(1\)
and are independent of the arrows. Hence their conditional avoidance
probability is
\[
 \exp\!\left(-\int_0^t|B_s|\dd s\right)
 =e^{-rt}\prod_{i=1}^r
   \exp\!\left(-\int_0^t X_i(s)\dd s\right).
\]
Define
\begin{equation}\label{eq:gap-functional}
 h_k(t)=\E\exp\!\left(-\int_0^t\min\{P_s,k\}\dd s\right)
 \quad(k\ge0).
\end{equation}
In particular, \(h_0(t)=1\). Taking expectations and using the
independence of the gap processes gives
\begin{equation}\label{eq:gap-factorization}
 \Pp(S\text{ is entirely uninfected at time }t)
 =e^{-rt}\prod_{i=1}^r h_{k_i}(t)
 =a^r\prod_{i=1}^r h_{k_i}(t).
\end{equation}

For the gap functional, fix \(k\ge1\), and let
\(E\sim\Exp(1)\) be the first arrival time of
\(P\). This is the first backward arrow that enlarges the ancestor
set within the gap. If \(E>t\), the integral in
\eqref{eq:gap-functional} is zero, contributing \(e^{-t}\) to
\(h_k(t)\). Conditional on \(E=s\le t\), the subsequent arrivals
form an independent rate-\(1\) Poisson process \(P'\), and
\[
 \int_0^t\min\{P_v,k\}\dd v
 =(t-s)+\int_0^{t-s}\min\{P'_u,k-1\}\dd u.
\]
The conditional expectation of the exponential in
\eqref{eq:gap-functional} is therefore
\(e^{-(t-s)}h_{k-1}(t-s)\). Since \(E\) has density \(e^{-s}\),
\begin{align}
 h_k(t)
 &=e^{-t}+\int_0^t e^{-s}e^{-(t-s)}h_{k-1}(t-s)\dd s\notag\\
 &=e^{-t}\left(1+\int_0^t h_{k-1}(u)\dd u\right).
 \label{eq:gap-recurrence}
\end{align}
Differentiating, and recalling that \(a=e^{-t}\), yields
\begin{equation}\label{eq:gap-differential-recurrence}
 h_k'(t)=-h_k(t)+a h_{k-1}(t),\qquad h_k(0)=1
 \quad(k\ge1),
\end{equation}
with the separate boundary condition \(h_0(t)\equiv1\).

For each vertex \(v\), let \(A_v\) be the event that \(v\) is
uninfected at time \(t\). Since complete infection means that none
of the events \(A_v\) occurs, inclusion--exclusion and
\eqref{eq:gap-factorization} give
\begin{align}
 F_m(t)
 &=1+\sum_{r=1}^m(-1)^r
       \sum_{\substack{S\subseteq V\\|S|=r}}
       \Pp\!\left(\bigcap_{v\in S}A_v\right)\notag\\
 &=1+\sum_{r=1}^m(-a)^r W_{m,r}(t),
 \label{eq:cycle-inclusion-exclusion}
\end{align}
where
\[
 W_{m,r}(t)=\sum_{\substack{S\subseteq V\\|S|=r}}
              \prod_{i=1}^r h_{k_i(S)}(t).
\]
The weight of \(S\) is well defined, since choosing a different
starting vertex only cyclically permutes the gap lengths.

To evaluate \(W_{m,r}\), count pairs \((S,v)\) with \(|S|=r\) and
\(v\in S\), assigning each pair the weight of \(S\). Each subset
has exactly \(r\) choices of the marked vertex, so their total
weight is \(rW_{m,r}(t)\). Alternatively, choose the marked vertex
in \(m\) ways and specify its ordered clockwise gaps
\((k_1,\ldots,k_r)\), where the entries are nonnegative and sum to
\(m-r\). These data determine a unique marked subset: starting at
the marked vertex, move successively \(k_1+1,\ldots,k_r+1\) steps
around the cycle. These increments are positive and sum to \(m\),
so the selected vertices are distinct and the last move returns to
the starting vertex. Thus
\begin{equation}\label{eq:marked-gap-count}
 rW_{m,r}(t)
 =m\!\sum_{\substack{k_1,\ldots,k_r\ge0\\
                     k_1+\cdots+k_r=m-r}}
       \prod_{i=1}^r h_{k_i}(t).
\end{equation}
This count also holds when the gap sequence is periodic: the
\(r\) marked vertices are distinct choices, even if some of the
resulting cyclic gap sequences coincide.
Put \(H(x,t)=\sum_{k\ge0}h_k(t)x^k\), and let \([x^j]\) denote
coefficient extraction. Expanding the product of \(r\) copies of
\(H\) gives
\[
 [x^{m-r}]H(x,t)^r
 =\sum_{\substack{k_1,\ldots,k_r\ge0\\
                  k_1+\cdots+k_r=m-r}}
       \prod_{i=1}^r h_{k_i}(t).
\]
Consequently, \(W_{m,r}(t)=(m/r)[x^{m-r}]H(x,t)^r\), and
\eqref{eq:cycle-inclusion-exclusion} becomes
\begin{equation}\label{eq:cycle-gap-cdf}
 F_m(t)=1+\sum_{r=1}^m\frac mr(-a)^r
                    [x^{m-r}]H(x,t)^r.
\end{equation}
For \(m=1\), its right-hand side is \(1-a\), since \(h_0=1\),
so the same identity holds with the stated convention for \(F_1\).

Divide \eqref{eq:cycle-gap-cdf} by \(m\), multiply by \(x^m\),
and sum over \(m\ge1\). Using \(z\) as the coefficient-extraction
variable, we obtain
\begin{align*}
 \sum_{m\ge1}\frac{F_m(t)}m x^m
 &=-\log(1-x)+\sum_{r\ge1}\frac{(-a)^r}{r}
       \sum_{m\ge r}x^m[z^{m-r}]H(z,t)^r\\
 &=-\log(1-x)+\sum_{r\ge1}\frac{(-axH(x,t))^r}{r}\\
 &=-\log(1-x)-\log(1+axH(x,t))\\
 &=-\log K(x,t),
 \qquad K(x,t)=(1-x)(1+axH(x,t)).
\end{align*}
All rearrangements are valid as formal power series, since each
coefficient receives only finitely many contributions.
The recurrence for \(h_k\) implies
\[
 H_t=1-(1-ax)H,\qquad
 K_t=2(1-x)-(2-ax)K,\qquad K(x,0)=1.
\]
Writing \(K=1-\sum_{j\ge1}c_jx^j\) yields
\[
 c_1=1-a,\qquad c_j'=-2c_j+ac_{j-1},\quad c_j(0)=0\quad(j\ge2).
\]
Equivalently, as functions of \(a\),
\(\frac{d}{da}(c_j/a^2)=-c_{j-1}/a^2\).
Repeated integration from \(a\) to \(1\), starting with
\(c_1/a^2=(1-a)/a^2\), proves \eqref{eq:cycle-coefficients}.
\end{proof}

For the remainder of the section, fix \(t>0\) and put
\[
 b=\sqrt{q(t)},\qquad
 p=\frac{1-e^{-t}}b,\qquad \delta=1-p,\qquad
 g_m=\frac{F_m(t)}{b^m}.
\]
The coefficient of \(x^2\) in
\eqref{eq:cycle-generating-function} gives
\(b^2=c_1^2+2c_2\). Since \(0<c_2\le c_1^2/2\) by
\eqref{eq:cycle-coefficients}, we have
\begin{equation}\label{eq:p-delta-range}
 \frac1{\sqrt2}\le p<1,\qquad 0<\delta<p.
\end{equation}

\begin{lemma}\label{lem:cycle-coefficient-comparison}
Define \(\Delta_j\) by
\[
 C_t(x/b)=1-(1-x)e^{\delta x}+\Delta(x).
\]
Then
\[
 \Delta(x)=\sum_{j\ge3}\Delta_jx^j,\qquad \Delta_j>0,
 \qquad
 \Delta_4\le\frac{\Delta_3}{2},\quad
 \Delta_5\le\frac{\Delta_3}{4}.
\]
\end{lemma}

\begin{proof}
Set \(\alpha=e^{-t}/b\). Substitution \(v=by\) in
\eqref{eq:cycle-coefficients} gives
\[
 \frac{c_j}{b^j}
 =\int_0^p\frac{y^{j-2}}{(j-2)!}\,w(y)\dd y,\qquad
 w(y)=\frac{\alpha^2(p-y)}{(\alpha+y)^2}\quad(j\ge2).
\]
Let \(\mu\) have density \(w(y)\) on \([0,p]\), and let \(\nu\)
have density \(p+y\) on \([0,\delta]\). Their masses agree:
\[
 \mu([0,p])=\frac{c_2}{b^2}
 =\frac{1-p^2}{2}=p\delta+\frac{\delta^2}{2}
 =\nu([0,\delta]).
\]
For \(0<y<\delta\), we have \(w(y)\le p-y<p+y\).
Thus \(\mu-\nu\) is negative below \(\delta\), positive above
\(\delta\), and has total mass zero. The integral of an increasing
function against \(\mu\) is therefore at least its integral against
\(\nu\), with strict inequality when the function is strictly
increasing. Applying this observation to the monomials gives
\[
 \frac{c_j}{b^j}\ge
 \int_0^\delta\frac{y^{j-2}}{(j-2)!}(p+y)\dd y
 =\frac{\delta^{j-1}(j-\delta)}{j!},
\]
with equality for \(j=2\) and strict inequality for \(j\ge3\).
The expression on the right is the coefficient of \(x^j\) in
\(1-(1-x)e^{\delta x}\), whose linear coefficient is \(p=c_1/b\).
This proves the assertions about the coefficients of \(\Delta\).

For the remaining bounds, put
\(L(y)=\nu([0,y])-\mu([0,y])\). Then \(L\ge0\), and integration
by parts gives
\[
 \Delta_j=\int_0^p\frac{y^{j-3}}{(j-3)!}L(y)\dd y
 \qquad(j\ge3).
\]
On \([\delta,p]\), the function \(L(y)=\mu([y,p])\) is
nonincreasing. Since \(\delta+p=1\) and \(\delta<1/2\), pairing
\(y\) with
\(1-y\) on \([\delta,p]\) gives
\[
\begin{aligned}
 \Delta_4-\frac12\Delta_3
 &=\int_0^\delta\left(y-\frac12\right)L(y)\dd y\\
 &\quad+\int_\delta^{1/2}\left(y-\frac12\right)
       \bigl(L(y)-L(1-y)\bigr)\dd y
 \le0.
\end{aligned}
\]
Both integrals are nonpositive: \(L\ge0\) on the first interval,
and \(L(y)\ge L(1-y)\) on the second. Finally, \(0\le y\le p<1\)
implies \(\Delta_5\le\Delta_4/2\le\Delta_3/4\).
\end{proof}

\begin{lemma}\label{lem:cycle-comparisons}
We have \(g_2=1\), \(g_m>1\) for \(m\ge3\), and
\[
 1<g_3<g_5<g_3^3,\qquad
 g_m>g_3\ (m\ge5),\qquad g_m>g_5\ (m\ge8).
\]
\end{lemma}

\begin{proof}
The coefficients
\[
 s_j=[x^j]\frac{e^{-\delta x}}{1-x}
     =\sum_{i=0}^j\frac{(-\delta)^i}{i!}
\]
lie between \(p=1-\delta\) and \(1\), by the alternating-series
bounds. Put
\[
 D(x)=\Delta(x)\frac{e^{-\delta x}}{1-x},
 \qquad E(x)=-\log(1-D(x)).
\]
These series have nonnegative coefficients and start in degree \(3\).
Since
\[
 1-C_t(x/b)=(1-x)e^{\delta x}(1-D(x)),
\]
the generating-function identity gives
\[
 g_m=1+m[x^m]E(x)\quad(m\ge2).
\]
This proves \(g_2=1\) and \(g_m>1\) for \(m\ge3\). In particular,
\[
 g_3-1=3\Delta_3,\qquad
 g_5-1=5(s_2\Delta_3+p\Delta_4+\Delta_5).
\]
The coefficientwise inequality \(E\ge D\) and
\eqref{eq:p-delta-range} now give
\[
 \begin{aligned}
 g_m-1&\ge mp\Delta_3>3\Delta_3=g_3-1
                         &&(m\ge5),\\
 g_m-1&\ge mp(\Delta_3+\Delta_4+\Delta_5)
          >5(\Delta_3+\Delta_4+\Delta_5)\ge g_5-1
                         &&(m\ge8).
 \end{aligned}
\]
Finally, Lemma~\ref{lem:cycle-coefficient-comparison} yields
\[
 g_5-1
 \le5\left(1+\frac12+\frac14\right)\Delta_3
 =\frac{35}{4}\Delta_3
 <9\Delta_3=3(g_3-1)<g_3^3-1.
\]
\end{proof}

\begin{proof}[Proof of Theorem~\ref{thm:main}]
Fix \(t>0\). By \eqref{eq:auxiliary-cycle-law}, the completion-time
distribution function of a cycle \(C_\ell\), divided by \(b^\ell\),
is \(g_\ell\) for
odd \(\ell\) and \(g_{\ell/2}^2\) for even \(\ell\).
Lemma~\ref{lem:cycle-comparisons} therefore permits replacing a cycle
by four-cycles together with, respectively, no additional component,
one \(C_5\), two \(C_3\)'s, or one \(C_3\), according as its length
is congruent to \(0,1,2,3\) modulo \(4\). Each replacement preserves
the order and does not increase the distribution function. The
inequality is strict unless the original
cycle is \(C_3,C_4,C_5\), or \(C_6\); replacing \(C_6\) by
\(2C_3\) leaves the distribution function unchanged.

The resulting graph is \(uC_3\sqcup wC_4\sqcup vC_5\), whose
normalized distribution function is
\(g_3^u g_5^v\). The inequalities \(1<g_3<g_5<g_3^3\) give three
further strict replacements:
\[
 3C_3\longrightarrow C_5\sqcup C_4,\qquad
 2C_5\longrightarrow2C_3\sqcup C_4,\qquad
 C_3\sqcup C_5\longrightarrow2C_4.
\]
Each replacement increases the number of four-cycles, so iteration
terminates. In the resulting graph, \(u\le2\), \(v\le1\), and
\(uv=0\). The congruence \(3u+5v\equiv n\pmod4\) determines
\[
 (u,v)=(0,0),(0,1),(2,0),(1,0)
\]
in the four residue classes, respectively. These give exactly the
extremal graphs stated in the theorem, allowing for the equality
\(C_6\leftrightarrow2C_3\).

By component independence, all the comparisons of distribution
functions can be multiplied. Every graph outside the stated equality
cases requires at least one strict replacement, so its distribution
function is strictly larger at every \(t>0\).
Integrating the tail inequalities and applying
Lemma~\ref{lem:poisson} proves the expectation statement and its
equality cases.
\end{proof}

\begin{corollary}\label{cor:degree-two-universal-cdf}
For every finite simple \(2\)-regular graph \(G\) on \(n\) vertices,
\[
 \Pp(T_G\le t)\ge q(t)^{n/2}\qquad(t>0).
\]
Equality at any positive \(t\) holds if and only if all components
of \(G\) are four-cycles.
\end{corollary}

\begin{proof}
Apply \(g_2=1\) and \(g_m>1\) for \(m\ge3\) to
\eqref{eq:auxiliary-cycle-law}, and multiply over the components.
\end{proof}

Thus every extremal graph has all but at most six vertices in
four-cycle components. Among orders divisible by three, triangle
unions are extremal exactly
at orders three and six. The smallest counterexample in degree two
therefore has nine vertices:
\begin{equation}\label{eq:nine-vertex}
 \cc(C_4\sqcup C_5)>\cc(3K_3).
\end{equation}
Theorem~\ref{thm:main} also gives counterexamples at every larger
order divisible by three. At orders divisible by twelve, the comparison
has the following explicit form.

\begin{corollary}\label{cor:counterexample}
For every integer \(k\geq1\),
\begin{equation}\label{eq:family}
 \cc(3kC_4)>\cc(4kK_3).
\end{equation}
Both graphs are simple and \(2\)-regular on \(12k\) vertices, and
\begin{equation}\label{eq:strict-cdf}
 \Pp(T_{3kC_4}\leq t)<\Pp(T_{4kK_3}\leq t)
 \qquad(t>0).
\end{equation}
In particular,
\begin{align}
 \cc(3C_4)&=\frac{390968681}{16200000}
             \approx24.134,\label{eq:exact-G}\\
 \cc(4K_3)&=\frac{222738867019}{9685021500}
             \approx22.998.\label{eq:exact-H}
\end{align}
\end{corollary}

The first exact value also equals the expected number of rolls of a
fair six-sided die needed to see every face twice; it is recorded in
\cite[Problem~9]{Conroy2018}.

\begin{proof}[Proof of Corollary~\ref{cor:counterexample}]
The strict comparisons follow from Theorem~\ref{thm:main}, since
\(3kC_4\) is the unique extremal graph at order \(12k\), whereas
\(4kK_3\) is not isomorphic to it. For \(k=1\), the tail-integral
identity \eqref{eq:tail} and Lemma~\ref{lem:poisson} give
\begin{align*}
 \cc(3C_4)&=12\int_0^\infty\bigl(1-q(t)^6\bigr)\dd t,\\
 \cc(4K_3)&=12\int_0^\infty\bigl(1-r(t)^4\bigr)\dd t.
\end{align*}
Expanding the integrands and applying
\[
 \int_0^\infty t^j e^{-at}\dd t=\frac{j!}{a^{j+1}}
 \qquad(a>0,\ j=0,1,2,\ldots)
\]
yields \eqref{eq:exact-G} and \eqref{eq:exact-H}.
\end{proof}

\begin{remark}
The distributional comparisons concern the Poissonized completion
times. Their expectation consequences for the discrete process follow
from \eqref{eq:poisson}.
\end{remark}

\section{Neighborhood structure and extremality}
\label{sec:neighborhood-tails}\label{sec:bounded-components}

For a nonempty vertex set \(U\), put
\(N_G(U)=\bigcup_{u\in U}N_G(u)\), using open neighborhoods.
Partition the vertices into equivalence classes according to equality
of their open neighborhoods.
If \(G\) has minimum degree \(\delta\geq1\), let
\(\tau_\delta(G)\) be the largest size of such a class consisting of
vertices of degree \(\delta\).

\begin{theorem}\label{thm:tail-classes}
Let \(G\) be a finite simple graph with minimum degree
\(\delta\geq1\), and put \(\tau=\tau_\delta(G)\).
There is a constant \(c_G>0\) such that
\[
 \Pp(T_G>t)\sim c_G t^{\tau-1}e^{-\delta t}
 \qquad(t\to\infty).
\]
\end{theorem}

\begin{proof}
When the uncovered set is \(U\), a ring is successful exactly when
its vertex belongs to \(N_G(U)\). The total rate of the next covering
and the rate at which a particular \(u\in U\) is covered are
\[
 \lambda(U)=|N_G(U)|\geq\delta,
 \qquad
 a_U(u)=\sum_{v\in N_G(u)}\frac{1}{|N_G(v)\cap U|}>0.
\]
In particular, every ordering of the vertices has positive probability
of occurring as the order in which they are covered.

Equality \(\lambda(U)=\delta\) holds precisely when \(U\) is contained
in a single class of degree-\(\delta\) vertices with the same open
neighborhood. Indeed, if \(|N_G(U)|=\delta\), then every \(u\in U\) satisfies
\(N_G(u)\subseteq N_G(U)\) and \(|N_G(u)|\geq\delta\), forcing
\(N_G(u)=N_G(U)\). The converse follows immediately.
Once a covering order reaches such a set, all subsequent nonempty
uncovered sets remain in the same class. Hence the number of holding
times of rate \(\delta\) along that order is at most \(\tau\).
This maximum is attained by first covering all vertices outside a
class of size \(\tau\).

Condition on a covering order \(\pi\). By
Lemma~\ref{lem:finite-state-chains}, the holding time at each state
is exponential with rate \(\lambda(U)\), independently of the next
transition. Thus the successive holding times, conditional on
\(\pi\), are independent. If the order contains \(k\geq1\)
holding times of rate \(\delta\), its conditional completion time is
\(X+Y\), where \(Y\) is a sum of \(k\) independent
\(\Exp(\delta)\) variables and \(X\) is an independent sum of
exponential variables of rates greater than \(\delta\). Consequently
\[
 \Pp(X+Y>t)
 \sim\frac{\delta^{k-1}}{(k-1)!}
       \E e^{\delta X}\,t^{k-1}e^{-\delta t}.
\]
To verify this asymptotic, condition on \(X\) and apply the Erlang
tail formula from Lemma~\ref{lem:tail-exponential}, namely
\(\Pp(Y>s)=e^{-\delta s}
\sum_{j=0}^{k-1}(\delta s)^j/j!\).
After division by \(t^{k-1}e^{-\delta t}\), the highest power
converges by dominated convergence with dominator
\(e^{\delta X}\), and the lower powers vanish.
All rates in \(X\) are integers at least \(\delta+1\), so
\(\E e^{(\delta+1/2)X}<\infty\); this also makes the contribution
from \(X>t\) negligible. An order containing no rate-\(\delta\)
holding time has an exponentially smaller tail by the same bound.
Finally, average over the finitely many covering orders. Their probabilities
are positive, the largest possible \(k\) is \(\tau\), and every
order attaining it has a positive leading coefficient. This proves
the result and the positivity of \(c_G\).
\end{proof}
For a \(d\)-regular graph write \(\tau(G)=\tau_d(G)\).

\begin{corollary}\label{cor:regular-tail-classes}
Let \(G\) be a finite simple \(d\)-regular graph, \(d\geq2\).
Then \(1\leq\tau(G)\leq d\), with equality \(\tau(G)=d\)
if and only if \(G\) has a \(K_{d,d}\) component.
If it has exactly \(b\geq1\) such components, then
\[
 \Pp(T_G>t)
 \sim\frac{2b\,d^{d-1}}{(d-1)!}\,t^{d-1}e^{-dt}.
\]
In particular, \(K_{d,d}\) is the unique connected \(d\)-regular
graph whose completion-time tail has a polynomial factor of degree
\(d-1\).
\end{corollary}

\begin{proof}
Let \(C\) be a class of vertices with the same open neighborhood
\(D\). Every vertex of \(D\) is adjacent to all of \(C\), so
\(|C|\leq d\). If \(|C|=d\), then \(|D|=d\),
\(C\cap D=\varnothing\), and the edges between the two classes
account for every edge incident with a vertex of \(C\cup D\).
Consequently, these vertices form a \(K_{d,d}\) component.
The converse is immediate.

In a \(K_{d,d}\) component, each class must generate \(d\) clock
rings to cover the opposite class. These two times are independent
sums of \(d\) variables of distribution \(\Exp(d)\). Thus
Lemma~\ref{lem:tail-exponential} gives, with
\[
 s_d(t)=e^{-dt}\sum_{j=0}^{d-1}\frac{(dt)^j}{j!},
\]
the component survival probability is
\[
 2s_d(t)-s_d(t)^2
 \sim\frac{2d^{d-1}}{(d-1)!}\,t^{d-1}e^{-dt}.
\]
Every other component has \(\tau\leq d-1\) and hence a lower-order
tail by Theorem~\ref{thm:tail-classes}. Independence of components
gives the stated coefficient.
\end{proof}

\begin{corollary}\label{cor:eventual-tail-extremality}
Fix \(d\geq2\) and a positive integer \(h\) divisible by \(2d\),
and put \(B=(h/(2d))K_{d,d}\). For all sufficiently large \(t\),
\(B\) uniquely maximizes \(\Pp(T_G>t)\), up to isomorphism,
among the \(d\)-regular graphs on \(h\) vertices.
\end{corollary}

\begin{proof}
If \(\tau(G)<d\), then
\(\Pp(T_G>t)/\Pp(T_B>t)\to0\). Otherwise the ratio tends to
\(2db/h\), where \(b\) is the number of \(K_{d,d}\) components
of \(G\). This limit equals \(1\) exactly when \(G\cong B\).
For every other graph the limit is strictly smaller than \(1\).
Since there are only finitely many graphs on \(h\) vertices, a common
threshold for \(t\) can be chosen.
\end{proof}

For regular graphs of bounded component order, the corresponding
extremal statement holds for the expected completion time. The
component multiplicities may vary arbitrarily with the total order.

\begin{theorem}\label{thm:bounded-component-extremality}
Fix integers \(d\geq2\) and \(L\geq2d\). For all sufficiently
large \(n\) divisible by \(2d\), the graph
\(B_n=(n/(2d))K_{d,d}\) uniquely maximizes \(\cc(G)\), up to
isomorphism, among all \(d\)-regular graphs on \(n\) vertices
whose components have at most \(L\) vertices.
More precisely, if \(s(G)\) counts the vertices lying outside
\(K_{d,d}\) components, then
\[
 \cc(B_n)-\cc(G)\geq\frac{e^{-d}}{32}\,s(G)
\]
for every sufficiently large admissible \(n\), uniformly over all
such \(G\).
\end{theorem}

\begin{proof}
Write \(F_H(t)=\Pp(T_H\leq t)\),
\(Q(t)=F_{K_{d,d}}(t)\), and \(S(t)=1-Q(t)\).
For each connected \(d\)-regular graph
\(H\not\cong K_{d,d}\) of order at most \(L\),
Corollary~\ref{cor:regular-tail-classes} gives
\((1-F_H(t))/S(t)\to0\). Since there are only finitely many such
component types, a common \(t_0>0\) can be chosen so that
\[
 1-F_H(t)\leq
 \min\left\{\frac12,\frac{|V(H)|}{8d}S(t)\right\}
 \qquad(t\geq t_0)
\]
for every such \(H\). Using
\(\log(1-x)\geq-2x\) for \(0\leq x\leq1/2\) and
\(-\log(1-y)\geq y\) for \(0\leq y<1\), we obtain
\[
 \log F_H(t)-\frac{|V(H)|}{2d}\log Q(t)
 \geq\frac{|V(H)|}{4d}S(t)
 \qquad(t\geq t_0).
\]

Put \(N=n/(2d)\) and \(s=s(G)\). Multiplying component
distribution functions and summing these inequalities gives
\[
 F_G(t)\geq Q(t)^N
       \exp\!\left(\frac{s}{4d}S(t)\right)
 \qquad(t\geq t_0).
\]
Since \(e^x-1\geq x\), the tail-integral identity
\eqref{eq:tail} implies
\begin{align*}
 \E T_{B_n}-\E T_G
 &=\int_0^\infty\bigl(F_G(t)-Q(t)^N\bigr)\dd t\\
 &\geq-t_0Q(t_0)^N
       +\frac{s}{4d}\int_{t_0}^\infty Q(t)^NS(t)\dd t.
\end{align*}
Here the possible negative contribution before \(t_0\) is bounded
using \(F_G\geq0\) and the monotonicity of \(Q\).

Choose \(u_N\) with \(S(u_N)=1/N\). Such a point exists for
\(N>1\), since \(S\) is continuous and strictly decreasing from
\(1\) to \(0\). We have \(u_N\to\infty\), and the explicit
tail of \(K_{d,d}\) gives
\(S(u_N+1)/S(u_N)\to e^{-d}\).
For all sufficiently large \(N\), the interval
\([u_N,u_N+1]\) lies above \(t_0\), and on this interval
\[
 S(t)\geq\frac{e^{-d}}{2N},
 \qquad
 Q(t)^N\geq(1-1/N)^N\geq\frac14.
\]
Consequently
\[
 \E T_{B_n}-\E T_G
 \geq\frac{s e^{-d}}{32dN}-t_0Q(t_0)^N.
\]
Multiplication by \(n=2dN\), using Lemma~\ref{lem:poisson}, yields
\[
 \cc(B_n)-\cc(G)
 \geq\frac{s e^{-d}}{16}-nt_0Q(t_0)^N.
\]
For \(G\not\cong B_n\), the integer \(s\) is positive and
divisible by \(2d\), so \(s\geq2d\).
Since \(0<Q(t_0)<1\), the last error tends to zero exponentially
in \(n\), independently of \(G\). For all sufficiently large
\(n\), it is at most \(s e^{-d}/32\), proving the asserted bound
and strict extremality. When \(G\cong B_n\), both sides are zero.
\end{proof}

\section{Clique unions and uniform bounds}\label{sec:clique-unions}

\begin{proposition}\label{prop:all-degrees}
Fix an integer \(d\geq2\). For \(n\) divisible by both \(2d\) and \(d+1\),
put
\[
 B_{n,d}=\frac{n}{2d}K_{d,d},
 \qquad
 H_{n,d}=\frac{n}{d+1}K_{d+1}.
\]
As \(n\to\infty\) through these common multiples,
\begin{equation}\label{eq:general-gap}
 \cc(B_{n,d})-\cc(H_{n,d})
 =\frac{(d-1)n}{d}\log \log  n+O_d(n).
\end{equation}
In particular, \(\cc(B_{n,d})>\cc(H_{n,d})\) for all sufficiently
large such \(n\).
\end{proposition}

\begin{proof}
Partition \(B_{n,d}\) into the two vertex classes of each of its
complete bipartite components. There are \(m=n/d\) classes, each of
size \(d\). A uniform vertex sample selects each class with probability
\(1/m\). A sample from one class covers a vertex of its opposite class
until the latter is full. The graph is therefore covered exactly when
each of the \(m\) classes has been sampled at least \(d\) times. Thus,
as discrete random variables,
\begin{equation}\label{eq:collector-identity}
 N_{B_{n,d}}\ \stackrel{\mathrm{d}}{=}\ D_{n/d,d}.
\end{equation}
Applying Theorem~\ref{thm:newman-shepp} gives
\begin{equation}\label{eq:bipartite-asymptotic}
 \cc(B_{n,d})
 =\frac nd\bigl(\log (n/d)+(d-1)\log \log (n/d)+O_d(1)\bigr).
\end{equation}

For a single \(K_{d+1}\), every ring is successful while at least two
vertices remain uncovered. Once only one uncovered vertex remains,
a ring is successful precisely when it occurs at one of its \(d\)
neighbors. The strong Markov property in
Lemma~\ref{lem:poisson-clocks}, applied at the time of the \(d\)th
covering, shows that the conditional distribution of the final waiting
time is \(\Exp(d)\), independently of the preceding history.
Consequently, its continuous completion time is \(S=X+Y\), where
\(X\) is the sum of \(d\) independent \(\Exp(d+1)\) variables and
\(Y\sim\Exp(d)\) is independent of \(X\). Hence
\begin{equation}\label{eq:clique-tail}
 e^{-dt}\leq\Pp(S>t)
 \leq e^{-dt}\E e^{dX}=(d+1)^d e^{-dt}
 \qquad(t\geq0).
\end{equation}
The lower bound follows from \(S\geq Y\). For the upper bound,
condition on \(X\) and use
\(\Pp(Y>t-X\mid X)\leq e^{-d(t-X)}\), which remains valid when
\(X>t\). For \(Z\sim\Exp(d+1)\),
\[
 \E e^{dZ}=(d+1)\int_0^\infty e^{-s}\dd s=d+1.
\]
Independence of the \(d\) summands in \(X\) therefore gives
\(\E e^{dX}=(d+1)^d\).
Let \(h=n/(d+1)\), and let \(M_h\) be the maximum of \(h\)
independent copies of \(S\). By component independence,
\(T_{H_{n,d}}\) has the same distribution as \(M_h\).
The lower bound in \eqref{eq:clique-tail} compares \(M_h\) with
the maximum of \(h\) independent \(\Exp(d)\) variables, whose
expectation is \(H_h/d\) by Lemma~\ref{lem:tail-exponential}, where
\(H_h=\sum_{j=1}^{h}1/j\). For the upper bound, apply the union
bound and integrate the resulting tail estimate to obtain
\[
 \frac{H_h}{d}\leq \E M_h
 \leq\int_0^\infty
       \min\{1,h(d+1)^d e^{-dt}\}\dd t
 =\frac{\log \bigl(h(d+1)^d\bigr)+1}{d}.
\]
Since \(H_h=\log h+O(1)\), Lemma~\ref{lem:poisson} yields
\begin{equation}\label{eq:clique-asymptotic}
 \cc(H_{n,d})=\frac nd\bigl(\log (n/(d+1))+O_d(1)\bigr).
\end{equation}
For fixed \(d\), the difference between \(\log(n/d)\) and
\(\log(n/(d+1))\) is the constant \(\log((d+1)/d)\), while
\(\log\log(n/d)=\log\log n+o(1)\). Subtracting
\eqref{eq:clique-asymptotic} from \eqref{eq:bipartite-asymptotic}
therefore proves \eqref{eq:general-gap}.
\end{proof}

The error term in \eqref{eq:general-gap} depends on the fixed degree
\(d\). In degree two, Theorem~\ref{thm:main} gives the exact
extremal classification at every order, and
Corollary~\ref{cor:counterexample} gives the strict comparison with
triangle unions already at order twelve. When the degree is allowed
to grow with the order, Alon
\cite[Section~6, Remark~3]{Alon2008} also suggested the uniform additive
estimate
\begin{equation}\label{eq:possible-bound}
 \cc(G)\leq n+O\bigl((n\log n)/d\bigr),
\end{equation}
with an absolute implied constant.

\begin{proposition}\label{prop:additive-obstruction}
For every integer \(m\geq2\),
\[
 0\leq
 \cc(2K_m)-\left(2m+
      \frac{2m}{4^m}\binom{2m}{m}\right)
 \leq\frac{3}{m-1}.
\]
In particular, \(\cc(2K_m)\geq2m+\sqrt m\). Consequently there
is no absolute constant \(C\) such that
\(\cc(G)\leq n+C(n\log n)/d\) for every \(d\)-regular graph
on \(n\) vertices.
\end{proposition}

\begin{proof}
Let \(\tau_m\) be the first step by which each component has
received at least \(m\) vertex samples. Every sample can cover at
most one vertex in its own component, so \(N_{2K_m}\geq\tau_m\).
After \(2m\) steps, the sample count in the first component is
\(B\sim\operatorname{Bin}(2m,1/2)\). The deficient component
still needs \(|B-m|\) samples, each arriving with probability
\(1/2\). Therefore
\[
 \E\tau_m=2m+2\E|B-m|
 =2m+\frac{2m}{4^m}\binom{2m}{m}.
\]
The last equality follows by symmetry and telescoping, using
\[
 (j-m)\binom{2m}{j}
 =m\left(\binom{2m-1}{j-1}-\binom{2m-1}{j}\right)
 \qquad(j>m).
\]
Writing \(a_m=4^{-m}\binom{2m}{m}\), induction gives
\(a_m\geq1/(2\sqrt m)\): the initial value is \(a_1=1/2\), and
\[
 \frac{a_{m+1}}{a_m}
 =\frac{2m+1}{2m+2}
 \geq\sqrt{\frac{m}{m+1}}.
\]
This proves the required lower bound. Since the graph has
\(n=2m\) and \(d=m-1\), its excess above \(n\) is at least
\(\sqrt m\), whereas \((n\log n)/d=O(\log m)\).

For the upper bound, construct the following coupling. Let
\(X_1,X_2\) be independent sums of
\(m-1\) independent variables of distribution \(\Exp(m)\), and let
\(E_1,E_2\sim\Exp(1)\) be independent of them and of each other.
The two clique completion times have the joint law of
\(S_i=X_i+E_i/(m-1)\), while the times of the \(m\)th clock
ring in the two components have the joint law of
\(Q_i=X_i+E_i/m\). Thus
\[
 0\leq\max_i S_i-\max_i Q_i
 \leq\frac{\max(E_1,E_2)}{m(m-1)}.
\]
Poissonization gives
\[
 \cc(2K_m)=2m\E\max_iS_i,
 \qquad \E\tau_m=2m\E\max_iQ_i.
\]
By Lemma~\ref{lem:tail-exponential},
\(\E\max(E_1,E_2)=3/2\), which proves the upper bound.
\end{proof}

The additive estimate therefore fails already for two disjoint
cliques because of the imbalance between the numbers of samples
received by the two components. Stirling's formula also gives
\[
 \cc(2K_m)=2m+2\sqrt{\frac m\pi}+O(m^{-1/2}).
\]
The weaker order estimate
\(\cc(G)=O(n+(n\log n)/d)\) is compatible with this family,
whose expected completion time is \(n+O(\sqrt n)=O(n)\).
 
\section*{Declaration of generative AI and AI-assisted technologies}
The authors used a large language model during the preparation of this manuscript
for assistance with mathematical exploration, checking arguments,
organization and language editing. The authors are responsible for
the mathematical statements, proofs, references and final text.

\bigskip
\begingroup
\small
\begin{samepage}
\noindent\textsc{School of Mathematical Science, Tianjin Normal
University, Tianjin 300387, P.R. China}\par
\noindent\hspace*{1em}\textit{E-mail address: }
\href{mailto:dengyuqiu1980@126.com}{\nolinkurl{dengyuqiu1980@126.com}}\par
\end{samepage}
\medskip
\begin{samepage}
\noindent\textsc{School of Mathematical Sciences, Tiangong University,
Tianjin 300387, P.R. China}\par
\noindent\hspace*{1em}\textit{E-mail address: }
\href{mailto:gqwang1979@aliyun.com}{\nolinkurl{gqwang1979@aliyun.com}}\par
\end{samepage}
\endgroup
\end{document}